\documentclass[twoside]{article}
\usepackage{amsthm,amsfonts,amsmath,amssymb}

\newcommand{\prob}{\mathbb P}
\newcommand{\ex}{\mathbb E\,}

\newcommand{\eq}{\begin{equation}}
\newcommand{\en}{\end{equation}}

\newcommand{\la}{\lambda}
\newcommand{\bfdelta}{\mbox{{\boldmath$\delta$}}}

\newcommand{\sign}{\operatorname{sign}}
\newcommand{\imi}{\mathtt{i}}

\newcommand{\intfloor}[1]{\left\lfloor#1\right\rfloor}
\newcommand{\intfloorsm}[1]{\lfloor#1\rfloor}
\newcommand{\fromto}{\mathord{\downarrow}}

\newcommand{\CM}{\mathcal{M}}

\newcommand{\CS}{\mathcal{S}}

\newcommand{\BN}{\mathbb{N}}
\newcommand{\BR}{\mathbb{R}}

\newtheorem{lemma}{Lemma}
\newtheorem{thm}[lemma]{Theorem}
\newtheorem{corollary}[lemma]{Corollary}
\theoremstyle{remark}
\newtheorem*{remark*}{Remark}
\newtheorem*{example*}{Example}

\title{On the number of collisions in
$\Lambda$-coalescents\thanks{Research supported
 by the Netherlands Organisation for Scientific Research (NWO).}}
\author{Alexander~Gnedin and Yuri~Yakubovich\thanks{Postal address:
 Department of Mathematics, Utrecht University,
 Postbus 80010, 3508 TA Utrecht, The Netherlands. E-mail addresses:
 \texttt{gnedin@math.uu.nl}, \texttt{yakubovich@math.uu.nl}.}}
\date{\empty}

\begin{document}
\maketitle

\begin{abstract}
\noindent
 We examine the total number of collisions $C_n$ in the
$\Lambda$-coalescent process which starts with $n$ particles.
A linear growth and a stable limit law for $C_n$ are shown under the
assumption of a power-like
behaviour of the measure $\Lambda$ near $0$ with
exponent $0<\alpha<1$.
\end{abstract}
\noindent
{\large Keywords:} $\Lambda$-coalescent, stable laws.
\vskip0.3cm
\noindent
{\large 2000 Mathematics Subject Classification:} 60J10, 60K15

\section{Introduction}
A system of particles undergoes a random Markovian evolution according to the rules of the Pitman--Sagitov $\Lambda$-coalescent \cite{Pitman99, Sagitov99} if the only 
possible type of interaction is a {\it collision}  affecting  two or more particles that merge together to form a single particle.
When the total number of particles is $b\geq 2$, a collision affecting some $2\leq j\leq b$ particles occurs at the probability rate  
\begin{equation}\label{eq:lanjasint}
\la_{b,j}=\binom{b}{j}\int_0^1 x^{j-2}(1-x)^{b-j}\Lambda(dx),
\end{equation}
where $\Lambda$ is a given finite measure on $[0,1]$.
Linear time change allows to rescale $\Lambda$ by its total mass, 
making it a probability
measure, which is always supposed below.
Two important special cases are Kingman's coalescent~\cite{Kingman82}
with $\Lambda$ a unit mass at $0$ (when only binary collisions are possible),
and the Bolthausen--Sznitman coalescent~\cite{Bolthausen-Sznitman} with
$\Lambda$ the Lebesgue measure on $[0,1]$.
See \cite{BBS1, BBS2, DongGP, Moehle3} for recent work on the
$\Lambda$-coalescents and further references.

\par A quantity of considerable interest is the number of collisions $C_n$ which occur as the system  progresses from the initial state with $n$ particles
to the terminal state with a single particle.
Representing the coalescent process by a genealogical tree,  $C_n$ can be also 
understood as the number of non-leave nodes.  Asymptotic
properties of $C_n$ are sensitive functions of the behaviour of  $\Lambda$ near $0$. In this paper
we explore the class of measures which satisfy
\begin{equation}\label{ass2}
\Lambda\left([0,x]\right)=A x^\alpha+O(x^{\alpha+\varsigma})
~~~ {\rm as~~}x\downarrow 0,~~  {\rm with~~} 0<\alpha<1~~ {\rm and~~} \varsigma>0.   
\end{equation}
Under this assumption we show that $C_n\sim (1-\alpha)n$ as $n\to\infty$ (Lemma~\ref{lem:lln}) 
and that the law of $C_n$ approaches a 
completely asymmetric stable distribution of index $2-\alpha$
(Theorem~\ref{thm:numcoll}).

\par The same question for the
Bolthausen--Sznitman coalescent has been addressed recently
in~\cite{DIMR1, DIMR2}. This  can be viewed as a limiting case
of~\eqref{ass2} with $\alpha=1$. However, the technique of
\cite{DIMR1, DIMR2} is based on the particular form
of~$\Lambda$ in that case hence cannot be applied to the general $\Lambda$
satisfying~\eqref{ass2} with $\alpha=1$.
\par  If $\Lambda$ is a beta$(\alpha,2-\alpha)$ distribution with parameter $0<\alpha<2$,    
a time-reversal of the coalescent  describes the genealogy of a 
continuous-state branching process \cite{Birkner}. This connection was
exploited recently  to study  a small-time behaviour 
of $\Lambda$-coalescents~\cite{BBS1, BBS2} in the beta case.

\par We develop here a more robust and straightforward approach  based
on analysis of the decreasing Markov chain ${\CM}_n$ counting 
the number of particles. The number of collisions $C_n$ is the number
of steps needed for $\CM_n$ to reach the absorbing state~1 from
state~$n$. 
In Kingman's case $\CM_n$ has unit decrements, 
but in general the decrements of ${\CM}_n$ are not stationary, which is a major source
of difficulties preventing 
direct application
of the classical renewal theorems for step distributions with  infinite
variance~\cite{Feller}.
To override this obstacle  we show that when (\ref{ass2}) holds,
in a certain range ${\CM}_n$ can be bounded from above and below by
processes with stationary decrements.  It allows to approximate $C_n$,
and it happens that these bounds can be made tough enough
to derive the limit theorem.
Our method may be of interest in a wider context of the pure death
processes. 

\par By Schweinsberg's result \cite{Schweinsberg2000} a coalescent satisfying (\ref{ass2}) comes down from the infinity,
hence the number of particles existing at a fixed time is uniformly bounded whichever $n$.
Therefore the asymptotics of the number of collisions that occur prior some fixed time is the same as that of $C_n$.

\section{Markov chain ${\CM}_n$}\label{sec:gencoal}

Let ${\CM}_n$ be the 
Markov chain whose time ticks at the
collision events and the state 
coincides with the number of remaining particles.
Since no two collisions occur simultaneously, 
the number of collisions $C_n$ in the $\Lambda$-coalescent starting with $n$ particles
is the number of steps the Markov chain ${\CM}_n$ needs to proceed from
the initial state $n$ to the terminal state $1$.
Note that the number of particles decreases by $j-1$ when
 a collision affects $j$ particles, hence the probability of transition
from $b$ particles  to $b-j+1$ is 
\begin{equation}\label{eq:distJn}
q_b(j):=\frac{\la_{b,j}}{\la_b}\,,\qquad 
 2\leq j\leq b\,,
\end{equation} 
where $\lambda_b$ is the total 
 collision rate of $b$ particles 
\begin{equation}\label{eq:lanasint}
\la_b=\sum_{j=2}^b \la_{b,j}
=\int_0^1\frac{1-(1-x)^b-bx(1-x)^{b-1}}{x^2}\,\Lambda(dx).
\end{equation}

\par It is convenient to introduce the sequence of moments 
\begin{equation}\label{eq:defnu}
\nu_b:=\int_0^1 (1-x)^b\Lambda(dx)\,,~~~~~b=0,1,\ldots
\end{equation}
In view of $\la_{b,2}=\binom{b}{2}\nu_{b-2}$ the
rates $\la_{b,2}$ ($b=2,3,\dots$) uniquely  determine the whole
array $\la_{b,j}$, as one can also  
conclude from the consistency
relation 
$$
(b+1)\la_{b,j}=(b+1-j)\la_{b+1,j}+(j+1)\la_{b+1,j+1},
$$
which is equivalent to the integral representation of rates (\ref{eq:lanjasint}), see \cite{Pitman99}.

Simple
computation shows that the rates  can be derived from $\nu_b$'s as
\begin{equation}\label{lanj}
\la_{b,j}=\binom{b}{j}\sum_{s=0}^{j-2}(-1)^{j-s}
  \binom{j-2}{s}\nu_{b-2-s},
\end{equation}
and, from  
$\la_{b+1}-\la_b=b\nu_{b-1}$, we have
\begin{equation}\label{lan}
\la_b=\sum_{i=1}^{b-1}i\nu_{i-1}.
\end{equation}
Since the second difference of
$\int_0^1\frac{bx+(1-x)^b-1}{x^2}\Lambda(dx)$
is $\nu_b$, it also follows that
\begin{equation}\label{eq:expjn}
\sum_{j=2}^b(j-1)\la_{b,j}=\int_0^1\frac{bx+(1-x)^b-1}{x^2}
\Lambda(dx)
=\sum_{i=1}^{b-1}(b-i)\nu_{i-1}.
\end{equation}

\par We shall denote $J_b$ a random variable with distribution 
$${\mathbb P}(J_b=j)=q_b(j),$$
so the first decrement of ${\cal M}_n$ is distributed as $J_n-1$, and its mean value is 
\begin{equation}\label{defexJn}
\ex\bigl[J_n-1\bigr]=\sum_{j=2}^n(j-1)q_n(j)
=\frac{\sum_{i=1}^{n-1}(n-i)\nu_{i-1}}{\sum_{i=1}^{n-1}i\nu_{i-1}}
=n\frac{\sum_{i=1}^{n-1}\nu_{i-1}}{\sum_{i=1}^{n-1}i\nu_{i-1}}-1.
\end{equation}

\par Let $g(n,b)$ denote the  \textit{Green kernel}  equal to the
probability that the Markov chain $\CM_n$ ever visits state $b$.  
We have $g(n,n)=1$, and $g(n,1)=1$ since $1$ is the 
absorbing state  reached in at most ~$n-1$ steps. 
Decomposition over the first jump shows that the Green kernel
satisfies the recursion
\begin{equation}\label{eq:recGreen}
g(n,b)=\sum_{j=2}^{n-b+1}q_n(j)g(n-j+1,b),~~~~~~~~~~n>b\ge1.
\end{equation}
The moments of $C_n$ can be readily expressed in terms of the Green kernel.

\begin{lemma}\label{lem:momCn}
The first two moments of the number of collisions in the
$\Lambda$-coalescent started with $n$ particles are
\begin{align}\label{eq:momCn1}
\ex\bigl[C_n\bigr]&{}=\sum_{b=2}^n g(n,b)\,,\\ \label{eq:momCn2}
\ex\bigl[C_n^2\bigr]&{}=\sum_{b=2}^n g(n,b)\left(1+2\sum_{j=2}^{
b-1}g(b,j)\right)\,.
\end{align}
\end{lemma}
\begin{proof}
Formula (\ref{eq:momCn1}) is obvious since the number of collisions is the number of sites $b>1$ visited by ${\cal M}_n$.
Still, it is instructive to derive (\ref{eq:momCn1})  from the first-step decomposition
\begin{equation}\label{eq:disteqcn}
C_n=_d 1+C_{n-J_n+1},
\end{equation}
where in the RHS $J_n$  has distribution
\eqref{eq:distJn} and is 
independent from $C_1,\ldots, C_n$, with  $C_1=0$.
Taking expectations on both sides of~\eqref{eq:disteqcn} we obtain
\begin{equation}\label{eq:momCn:t1}
\ex\bigl[C_n\bigr]=h_n+\sum_{j=2}^nq_n(j)\ex\bigl[C_{n-j+1}\bigr],
\end{equation}
where $h_n\equiv 1$.
Replace now repeatedly $\ex\bigl[C_{n-1}\bigr]$,
$\ex\bigl[C_{n-2}\bigr]$, etc.\ using this 
recursion. Collecting coefficients at $h_b$ we see that it is the sum
over all decreasing paths from $n$ to~$b$ of probabilities of the path,
that is $g(n,b)$ by definition. Using relations $g(n,n)=1$
and $\ex\bigl[C_1\bigr]=0$ we arrive at~\eqref{eq:momCn1}.

\par
To calculate the second moment the equation~\eqref{eq:disteqcn} is squared, so
\[
C^2_n=_d 1+2 C_{n-J_n+1} +\left(C_{n-J_n+1}\right)^2,
\]
from which 
\begin{equation}\label{eq:momCn:t2}
\ex\bigl[C_n^2\bigr]= 1+2\sum_{j=2}^n q_n(j)\ex\bigl[C_{n-j+1}\bigr]
+\sum_{j=2}^n q_n(j)\ex\bigl[C_{n-j+1}^2\bigr]\,.
\end{equation}
This has the same structure as \eqref{eq:momCn:t1} with 
$h_n=1+2\sum_{j=2}^n q_n(j)\ex\bigl[C_{n-j+1}\bigr]$.  Expressing
$\ex\bigl[C_{n-j+1}\bigr]$ from \eqref{eq:momCn1} 
and using recursion~\eqref{eq:recGreen} yield
\[
h_n=1+2\sum_{j=2}^n q_n(j)
\sum_{b=2}^{n-j+1}g(n-j+1,b)
=1+2\sum_{b=2}^{n-1}\sum_{j=2}^{n-b+1}q_n(j)g(n-j+1,b)
=1+2\sum_{b=2}^{n-1}g(n,b)\,.
\]
Following the same line as above we get
\eqref{eq:momCn2}.
\end{proof}

\section{Asymptotics of the moments}\label{sec:asympt}

From now on we only consider   measures $\Lambda$ satisfying
\eqref{ass2}.
Standard Tauberian arguments show that in this case 
\begin{equation}\label{eq:estnu2}
\nu_n = A \Gamma(\alpha+1) n^{-\alpha} + O(n^{-\alpha-\varsigma'}),
\qquad\qquad n\to\infty\,.
\end{equation}
Here and henceforth
\[
\varsigma'=\min\{1,\varsigma\}\,.
\]
This behaviour will imply that the transition probabilities $q_n(j)$ stabilise
as $n\to\infty$ for each fixed $j$. 
The relevant asymptotics of $\la_n$ and
$\la_{n,j}$ appeared in~\cite[Lemma~4]{Bertoin-LeGall} under a less
restrictive assumption of regular variation,  but we need to explicitly
control the error term.

\begin{lemma}\label{lem:latail}
Suppose $\Lambda$ satisfies \eqref{ass2}. Then for $n$ sufficiently
large 
\[
\left|\sum_{j=m}^n \la_{n,j} -\frac{A\alpha}{2-\alpha}\, 
\frac{\Gamma(m+\alpha-2)}{\Gamma(m)}\, n^{2-\alpha}\right|
< c\,
\frac{\Gamma(m+\alpha+\varsigma'-2)}
{\Gamma(m)}\,n^{2-\alpha-\varsigma'}\,,
\]
uniformly in $m=2,\dots,n$. 
\end{lemma}


\begin{proof}
Introduce the truncated moment
\[
G_{-2}(x)=\int_x^1\frac{\Lambda(dy)}{y^2}\,.
\]
Integrating by parts we derive from~\eqref{ass2}  that for 
$x\to0$
\begin{align*}
G_{-2}(x)&{}=
\frac{A\alpha}{2-\alpha}\,x^{\alpha-2}
+O(\max\{x^{\alpha+\varsigma-2},1\})
\end{align*}
Rewriting  \eqref{eq:lanjasint} in terms of $G_{-2}$ 
and integrating by parts we obtain
\[
\sum_{j=m}^n \la_{n,j}
=-\int_0^1\sum_{j=m}^n\binom{n}{j}x^{j}(1-x)^{n-j}dG_{-2}(x)
=m\binom{n}{m}\int_0^1x^{m-1}(1-x)^{n-m}G_{-2}(x)\,dx\,,
\]
because the sum telescopes and the integrated terms vanish. 
Plugging the expansion of $G_{-2}$ gives
\[
\left|\sum_{j=m}^n \la_{n,j}
-\frac{A\alpha}{2-\alpha}\,\frac{\Gamma(m+\alpha-2)\Gamma(n+1)}{
  \Gamma(m)\Gamma(n+\alpha-1)} \right|<
c\,
\frac{\Gamma(m+\alpha+\varsigma-2)\Gamma(n+1)}{
\Gamma(m)\Gamma(n+\alpha+\varsigma-1) }
\]
for $\varsigma< 2-\alpha$, in which case 
the result follows from the familiar asymptotics of the gamma function
$\Gamma(n+\beta)/\Gamma(n)=n^\beta+O(n^{\beta-1})$ ($n\to\infty$). If
$\varsigma>1$ the error term in this expansion constitutes the main
part of error, yielding appearance of $\varsigma'$ instead of
$\varsigma$.
The case $\varsigma\geq 2-\alpha$ is treated in the same way.
\end{proof}

\begin{corollary}\label{cor:lan}
If measure $\Lambda$  satisfies~\eqref{ass2} then as $n\to\infty$
\begin{eqnarray}\label{obiglan}
\nonumber
\la_n&=&\frac{A\Gamma(\alpha+1)}{2-\alpha}  \,n^{2-\alpha}
  +O\left(n^{2-\alpha-\varsigma'}\right),\\ 
\la_{n,j}&=&\frac{A\alpha\Gamma(j+\alpha-2)}{j!}\,n^{2-\alpha}
  +O\left(n^{2-\alpha-\varsigma'}\right),\\
\nonumber
q_n(j)&=&(2-\alpha)\frac{(\alpha)_{j-2}}{j!}
+O\left(n^{-\varsigma'}\right)\,
\end{eqnarray}
 for every fixed $j$.   
\end{corollary}

\begin{proof} The formula for $\la_n$ follows by the direct application of
Lemma~\ref{lem:latail} with $m=2$.  Expression for $\la_{n,j}$
is a difference between two subsequent tail sums.  The ratio of these
quantities gives $q_n(j)$.
\end{proof}

Thus $J_n$ converge in distribution. The convergence in mean is also true.
Note that the mean of the limiting distribution of jumps $J_b-1$ is
\begin{equation}\label{eq:expoflim}
\sum_{j=2}^\infty
(j-1)(2-\alpha)\frac{(\alpha)_{j-2}}{j!}=\frac{1}{1-\alpha}\,.
\end{equation}

\begin{lemma}\label{lem:exJn}
If\/~\eqref{ass2} holds then the mean decrease of the number of
particles after collision satisfies
\[
\ex\bigl[J_n-1\bigr]=\frac{1}{1-\alpha}+O\left(n^{-\min\{1-\alpha,\,
\varsigma\}} \right).
\]
\end{lemma}

\begin{proof}
By assumption~\eqref{ass2} relation~\eqref{eq:estnu2} implies
the existence of constants $n_0, c$ such that
\[
\left|\nu_{n-1}-A\Gamma(\alpha+1) n^{-\alpha}\right|
<c n^{-\alpha-\varsigma'},\qquad\qquad n\ge n_0\,.
\]
Approximating sums by integrals yields, as $n\to\infty$, 
\begin{align*}
\sum_{i=1}^{n-1}\nu_{i-1}
&{}=\sum_{i=n_0}^{n-1}A\Gamma(\alpha+1)i^{-\alpha}
+O\left(\sum_{i=n_0}^{n-1}i^{-\alpha-\varsigma'}\right)
+\sum_{i=1}^{n_0-1}\nu_{i-1}\\
&{}=A\Gamma(\alpha+1)n^{1-\alpha}(1+O(1/n))\int_{n_0/n}^1x^{-\alpha}dx
+O(n^{1-\alpha-\varsigma'})+O(1)\\
&{}=\frac{A\Gamma(\alpha+1)}{1-\alpha}n^{1-\alpha}
+O(\max\{1,n^{1-\alpha-\varsigma'}\})=\frac{A\Gamma(\alpha+1)}{1-\alpha}
n^{1-\alpha}
+O(\max\{1,n^{1-\alpha-\varsigma}\})
\end{align*}
by definition of $\varsigma'$.
Substitution of this expression into~\eqref{defexJn}
and applying Corollary~\ref{cor:lan} finishes the proof.
\end{proof}

\smallskip 
\noindent
{\bf Example.}
It is possible to choose measure $\Lambda$ so that 
the decrement probabilities for $j<n$ are exactly the same as for the limiting distribution truncated at $n$,
in which case the envisaged limit theorem for $C_n$ follows readily from \cite{Feller}.
To achieve
\[
q_n(j)=(2-\alpha)\frac{(\alpha)_{j-2}}{j!}~~~
(j=2,\dots,n-1),~~~
q_n(n)=\sum_{j=n}^\infty (2-\alpha)\frac{(\alpha)_{j-2}}{j!}
=\frac{\Gamma(n+\alpha-1)}{n!\,\Gamma(\alpha)}\,
\]
one should take the measure
$$\Lambda(dx) = \alpha\left(1-\frac{\alpha}{2}\right)
x^{\alpha-1}dx +
\frac{\alpha}{2}\,\bfdelta_1(dx),$$ 
which is a mixture of
 beta$(\alpha,1)$
and a Dirac mass at~1.  Adding
$\bfdelta_1$ does not affect $\la_{n,j}$ for $j<n$, so the integration
in~\eqref{eq:lanjasint} yields
\[
\la_{n,j}
=\binom{n}{j}
\left(1-\frac{\alpha}{2}\right)\frac{\Gamma(j+\alpha-2)(n-j)!}
{\Gamma(n+\alpha-1)} ~~~~ (j=2,\dots,n-1),~~~~
\la_{n,n}
=\left(1-\frac{\alpha}{2}\right)\frac{\alpha}{n+\alpha-2}
+\frac{\alpha}{2}\,.
\]
Summation (or direct integration of~\eqref{eq:lanasint}) implies
the desired expression for $q_n(j)$.

\par 
That a positive mass at $1$ does not affect the asymptotics of $C_n$ is seen e.g. by observing that 
the probability of total collision implied by this mass is of the order smaller than $n^{-1}$, namely $q_n(n)=O(n^{\alpha-2})$. 
On the continuous time scale of the coalescent, the mass at $1$ is responsible for the total coalescence time
(independent of $n$), hence the insensibility of the asymptotics to $\Lambda(\{1\})$ may be explained by the effect of coming 
down from the infinity, as mentioned in Introduction. 

\par The example also demonstrates that taking minimum in the error term of
Lemma~\ref{lem:exJn} is necessary.  Indeed, $\varsigma'=1$, 
however direct calculation using~\eqref{eq:expoflim} shows that
\[
\ex\bigl[J_n-1\bigr]=\frac{1}{1-\alpha}-\sum_{j=n}^\infty(j-1)(2-\alpha)\frac{(\alpha)_{j-2}}{j!}
+n\,q_n(n)
=\frac{1}{1-\alpha}-\frac{n^{\alpha-1}}{(1-\alpha)\Gamma(\alpha)}(1+O(1/n))\,.
\]
So the error term is $O(n^{\alpha-1})$, and not $O(n^{-1})$.
\vskip0.5cm

\begin{lemma}\label{lem:Green} If\/~\eqref{ass2} holds
then there exists $\upsilon\in{]0,1[}$ such that\/ 
$\lim g(n,k)=1-\alpha$ as $n,k\to\infty$ in such a way that\/  $k\ge
n^\upsilon$ and $n-k\to\infty$.
\end{lemma}
The heuristics for this is that the jumps $J_b-1$ become almost
identically distributed for large $b$ and their mean is close to
$1/(1-\alpha)$. If the distributions were indeed the same  
with
mean~$1/(1-\alpha)$ then the Lemma would follow from the renewal theorem.
We postpone a rigorous proof to
Section~\ref{sec:bounds}.

\begin{lemma}\label{lem:lln}
If\/~\eqref{ass2} holds 
then
\[
C_n\sim(1-\alpha)n ~~~~(n\to\infty)
\]
in probability.
\end{lemma}

\begin{proof} 
The argument is based on formulas of Lemma~\ref{lem:momCn}.
Indeed, immediately from Lemma~\ref{lem:momCn}
and Lemma~\ref{lem:Green},  $\ex\bigl[C_n\bigr]\sim (1-\alpha)n$ as $n\to\infty$. 
Similarly, as $b\to\infty$, 
$\sum_{j=2}^{b-1}g(b,j)\sim (1-\alpha)b$ 
hence
\[
\ex\bigl[C_n^2\bigr]\sim \sum_{b=2}^n2g(n,b)(1-\alpha)b
\sim (1-\alpha)^2n^2.
\]
 Application of  Chebyshev's
inequality completes the proof. 
\end{proof}

\section{Stochastic bounds on the jumps}\label{sec:bounds}

In this section we construct stochastic bounds on the decrements 
$J_b-1$ 
of Markov chain $\CM_n$ in a range 
 $b=k,\dots,n$, 
to control  the asymptotic behaviour of $C_n$. 
Specifically, we 
find random variables $J_{n\fromto k}^+$ and
$J_{n\fromto k}^-$ to secure
 the distributional bounds  
\begin{equation}\label{eq:ineqJJJ}
J_{n\fromto k}^+\leq_d J_b\leq_d J_{n\fromto k}^-\,.
\end{equation}
uniformly in some range $b=k,\dots,n$.
 Here $\leq_d$ denotes the stochastic order, meaning that  two random variables
$X$ and $Y$ satisfy $X\leq_d Y$ iff $\prob[X\le t]\ge \prob[Y\le t]$ for
all~$t$.

Our approach to establishing the limit theorem for the number 
of collisions is based on constructing random variables $J^+_{n\fromto k}$
and $J^-_{n\fromto k}$ which on the one hand comply with~\eqref{eq:ineqJJJ}
and on the other hand yield the same limit distribution 
of the sum of their independent copies. These two requirements point in opposite directions, 
forcing an adequate choice of these random variables to be a
compromise.  
We define the distributions which
depend on parameters  $\gamma,\beta\in{]0,1[}$ and 
$\theta\in{]\beta,1[}$. The calibration of these constants will be done later. 
For  $n\ge k>0$ define
\begin{equation}\label{eq:defqpm}
\begin{split}
q^-_{n\fromto k}(j)&:=\begin{cases}
\frac{\la_{n,2}-n^{-\gamma}\la_n\left(3:\intfloorsm{n^\beta}
\right)
 +\la_n\left(\intfloorsm{n^\beta}+1:n\right)}{\la_n}
-2\!\!\max\limits_{\ell\in\{k,\dots,n\}}\frac{\la_\ell\left(\intfloorsm{
n^\beta}+1:\ell\right)}{\la_\ell},\quad
  & j=2,\\
\frac{\la_{n,j}(1+n^{-\gamma})}{\la_n}, 
  & j=3,\dots,\intfloor{n^\beta},
\!\!\!\!\!\!\!\!\!\!\\
2\!\!\max\limits_{\ell\in\{k,\dots,n\}}\frac{\la_\ell\left(\intfloorsm{
n^\beta}+1:\ell\right)}{\la_\ell}
-2\!\!\max\limits_{\ell\in\{k,\dots,n\}}\frac{\la_\ell\left(\intfloorsm{
n^\theta}+1:\ell\right)}{\la_\ell},\quad
  & j=\intfloor{n^\theta}\,,\\
2\!\!\max\limits_{\ell\in\{k,\dots,n\}}\frac{\la_\ell\left(\intfloorsm{
n^\theta}+1:\ell\right)}{\la_\ell},\quad
  & j=n,\\
0,&\text{otherwise},
\end{cases}\\
q^+_{n\fromto k}(j)&:=\begin{cases}
\frac{\la_{k,2}+n^{-\gamma}\la_k\left(3:\intfloorsm{n^\beta}
\right)+\la_k\left(\intfloorsm{n^\beta}+1:k\right)}
{\la_k},\qquad\qquad &j=2,\\
\frac{\la_{k,j}(1-n^{-\gamma})}{\la_k},&j=3,\dots,\intfloor{
n^\beta},\\
0,&\text{otherwise},
\end{cases}
\end{split}
\end{equation}
where 
\[
\la_n(m:k)=\sum_{j=m}^k\la_{n,j}\,.
\]
Note that $\sum_j q^+_{n\fromto k}(j)=\sum_jq^-_{n\fromto k}(j)=1$.
Moreover, $q^\pm_{n\fromto k}(j)$ are nonnegative for
large enough $n$ and $k$. Indeed, the inequality $q^+_{n\fromto k}(j)\ge0$
is obvious. Lemma~\ref{lem:latail} implies that if $n$ and $k$ are large
enough and $k>n^\beta$ then
\begin{equation}\label{eq:laelltail}
\frac{c_1}{n^{\beta(2-\alpha)}}\le 
\frac{\la_\ell\left(\intfloorsm{n^\beta}+1:\ell\right)}{\la_\ell}
\le \frac{c_2}{n^{\beta(2-\alpha)}}
\end{equation}
for some $c_2>c_1>0$ uniformly in $\ell\in
\left\{k,\dots,n\right\}$. Hence (\ref{eq:laelltail})
 holds for the maximum over these $\ell$, and   
it follows that $q^-_{n\fromto k}(j)\ge 0$.

Hence, quantities $q^\pm_{n\fromto k}(j)$ define some
probability distributions on $\BN$, at least for large enough $n$
and $k$. Let
$J_{n\fromto k}^+$ and $J_{n\fromto k}^-$ be random variables 
with these distributions, so
\begin{equation}\label{eq:defJpJm}
\prob[J_{n\fromto k}^+=j]=q^+_{n\fromto k}(j)
\qquad\text{and}\qquad
\prob[J_{n\fromto k}^-=j]=q^-_{n\fromto k}(j).
\end{equation}
\begin{lemma}\label{lem:bounds} 
Suppose that $\beta$, $\gamma$, $\theta$ 
and $\upsilon$ satisfy the inequalities
\begin{equation}\label{eq:ineqboundsexist}
1>\upsilon>\theta>\beta>\gamma/(2-\alpha)>0
\qquad\text{and}\qquad
\gamma<\frac{(\upsilon-\beta)(2-\alpha)\varsigma'}
{2-\alpha-\varsigma'}\,.
\end{equation}
Then the stochastic bounds \eqref{eq:ineqJJJ} hold for $n$ large enough
and $b,k$ in the range $\intfloor{n^\upsilon}\leq k\leq b\leq n$. 
\end{lemma}

\begin{proof}
By definition of the stochastic order, we need to show that for
all~$m$
\begin{equation}\label{eq:ineqJJJJ}
\prob\left[J_{n\fromto k}^+\ge m\right] \le 
\prob\left[J_b\ge m\right] \le 
\prob\left[J_{n\fromto k}^-\ge m\right]\,.
\end{equation}

The first inequality (\ref{eq:ineqJJJJ}) is clearly true for $m\ge
\intfloor{n^\beta}+1$
because the left-hand 
side is zero and for $m\le 2$ because both sides are~1. 
Suppose $3\le m\le n^\beta$, then the first inequality reads as
\begin{equation}\label{eq:bounds:t1}
\la_b\la_k\left(m:\intfloorsm{n^\beta}\right)
\left(1-n^{-\gamma}\right)
\le \la_k\la_b(m:b)\,.
\end{equation}
Since $b\ge k\ge \intfloor{n^\upsilon}$, taking $n$ sufficiently large
enables us to apply Lemma~\ref{lem:latail} and Corollary~\ref{cor:lan} to
get asymptotic
estimates  valid for all $b$ in the range $k\le b\le n$. From the definition of
$\la_k(m:\intfloor{n^\beta})$, Lemma~\ref{lem:latail} and the inequality
$$\frac{\Gamma(m+\alpha+\varsigma'-2)}{\Gamma(m)}\ge 
\frac{\Gamma(\intfloorsm{n^\beta}+\alpha+\varsigma'-2)}
{\Gamma(\intfloorsm {n^\beta})}~~~{\rm for ~~}m\le\intfloor{n^\beta}$$
we obtain
\begin{multline*}
\la_k\left(m:\intfloor{n^\beta}\right)
{}=\la_k\left(m :k\right)
-\la_k\left(\intfloor{n^\beta}+1: k\right)\\
{}=\frac{A\alpha}{2-\alpha}\left(\frac{\Gamma(m+\alpha-2)}
{\Gamma(m)}
 -\frac{\Gamma(\intfloorsm{n^\beta}+\alpha-1)}
   {\Gamma(\intfloorsm{n^\beta}+1)} \right)k^{2-\alpha}
+O\left(\frac{\Gamma(m+\alpha+\varsigma'-2)}{\Gamma(m)}
 k^{2-\alpha-\varsigma'} \right).
\end{multline*}
Hence we
rewrite the inequality as
\begin{align*}
&\left(\frac{A\alpha}{2-\alpha}\left(\frac{\Gamma(m+\alpha-2)}
{\Gamma(m)}
 -\frac{\Gamma(\intfloorsm{n^\beta}+\alpha-1)}
   {\Gamma(\intfloorsm{n^\beta}+1)} \right)k^{2-\alpha}
+O\left(\frac{\Gamma(m+\alpha+\varsigma'-2)}{\Gamma(m)}
 k^{2-\alpha-\varsigma'} \right)\right)\\
&\qquad\qquad\qquad\qquad
\times\left(\frac{A\Gamma(\alpha+1)}{2-\alpha}b^{2-\alpha} 
+O\left(b^{2-\alpha-\varsigma'}\right)\right) 
\left(1-n^{-\gamma}\right)\\
&\qquad\qquad\le
\left(\frac{A\Gamma(\alpha+1)}{2-\alpha}k^{2-\alpha} 
+O\left(k^{2-\alpha-\varsigma'}\right)\right) \\
&\qquad\qquad\qquad\qquad
\times\left(\frac{A\alpha}{2-\alpha}\,
\frac{\Gamma(m+\alpha-2)}{\Gamma(m)}
 b^{2-\alpha}+O\left(\frac{\Gamma(m+\alpha+\varsigma'-2)}{\Gamma(m)}
 b^{2-\alpha-\varsigma'} \right)\right).
\end{align*}
The leading terms on both sides cancel, and simplifying 
this inequality we are reduced to checking
\begin{multline*}
O\left(\frac{\Gamma(m+\alpha+\varsigma'-2)}{\Gamma(m)}k^{
2-\alpha-\varsigma' } b^ { 2-\alpha }\right)
+O\left(\frac{\Gamma(m+\alpha-2)}{\Gamma(m)}
b^{2-\alpha-\varsigma'}k^{2-\alpha }\right)
\\
\le
\frac{A^2\alpha\Gamma(\alpha+1)}{(2-\alpha)^2}\,
\left(\frac{\Gamma(\intfloorsm{n^\beta}+\alpha-1)}
{\Gamma(\intfloorsm{n^\beta}+1)}
b^{2-\alpha}k^{2-\alpha}
+\frac{\Gamma(m+\alpha-2)}{\Gamma(m)}b^{2-\alpha}k^{2-\alpha}n^{-\gamma}
\right).
\end{multline*}
(In the above lines we neglected certain lower order terms
using inequalities like $b^{2-\alpha}k^{2-\alpha-\varsigma'}\ge
b^{2-\alpha-\varsigma'}k^{2-\alpha}$.)
The desired inequality follows from the following two inequalities
\begin{align} \label{eq:bounds:t2}
c_1&{}\le 
\frac{c_3\Gamma(m)}{n^{\beta(2-\alpha)}\Gamma(m+\alpha+\varsigma'-2)}
k^{\varsigma'}
+\frac{\Gamma(m+\alpha-2)}{\Gamma(m+\alpha+\varsigma'-2)}
k^{\varsigma'}n^{-\gamma}\,,\\
\label{eq:bounds:t3}
c_2&{}\le 
\frac{c_3\Gamma(m)}{n^{\beta(2-\alpha)}\Gamma(m+\alpha-2)}
k^{\varsigma'}
+k^{\varsigma'}n^{-\gamma}\,,
\end{align}
with sufficiently large constants $c_1,c_2>0$ (twice the ratio of a constant
implied by the corresponding $O(\cdot)$ and the constant in the
right-hand side is enough) and $c_3\in{]0,1[}$.

The right-hand side of~\eqref{eq:bounds:t2} considered as a function of 
$m\in{[3,\intfloorsm{n^\beta}]}$ has a unique minimum which is attained
at $m'\sim
\left(\displaystyle\frac{\varsigma'}{c_3(2-\alpha-\varsigma')}\right)^{
1/(2-\alpha)}
n^{\beta-\gamma/(2-\alpha)}$ and has the value asymptotic to 
\[
\left(\frac{c_3\Gamma(m')}{n^{\beta(2-\alpha)}
\Gamma(m'+\alpha+\varsigma'-2) }
+\frac{\Gamma(m'+\alpha-2)}{\Gamma(m'+\alpha+\varsigma'-2)}n^{-\gamma}
\right)
k^{\varsigma'}
\sim c_4 n^{-\beta\varsigma'-\gamma(2-\alpha-\varsigma')/(2-\alpha)}
k^{\varsigma'}.
\]
Since $k\ge \intfloor{n^\upsilon}$ the RHS grows to infinity
once~\eqref{eq:ineqboundsexist} holds.

\par
In inequality~\eqref{eq:bounds:t3} we neglect the first summand
in the RHS and still have the function
$b^{1-\alpha}n^{-\gamma}\ge 
n^{\upsilon\varsigma'-\gamma}$ which grows to infinity with
$n$ once~\eqref{eq:ineqboundsexist} holds.  Thus the first inequality 
in~\eqref{eq:ineqJJJJ} holds for all sufficiently large~$n$.

\smallskip
The second inequality in~\eqref{eq:ineqJJJJ} is obvious for 
$m\ge\intfloor{n^\beta}+1$ and for $m\le 2$. Suppose $3\le m\le
n^\beta$.
The inequality can be rewritten as 
\begin{align*}
\frac{\la_b(m:b)}{\la_b}
\le{}&\frac{\la_n(m:n)}{\la_n}\left(1+n^{-\gamma}\right)
+\left(1-n^{-\gamma}\right) \max_{\ell\in\{k,\dots,n\}}   
  \frac{\la_\ell(\intfloorsm{n^\beta}+1:\ell)}{\la_\ell}\\
&\qquad+\left(1+n^{-\gamma}\right)\left( \max_{\ell\in\{k,\dots,n\}}   
  \frac{\la_\ell(\intfloorsm{n^\beta}+1:\ell)}{\la_\ell}
  -\frac{\la_n(\intfloorsm{n^\beta}+1:n)}{\la_n}\right)\,.
\end{align*}
The latter follows from a simpler inequality
\begin{equation}\label{eq:bounds:t5}
\la_n\la_b(m:b)\le \la_b\la_n(m:n)\left(1+n^{-\gamma}\right)
+\la_b\la_n\left(1-n^{-\gamma}\right)
\max_{\ell\in\{k,\dots,n\}}\frac{\la_\ell(\intfloorsm{n^\beta}+1:\ell)}{
\la_\ell}\,.
\end{equation}
Since $k>\intfloor{n^\beta}$, application of~\eqref{eq:laelltail}
implies 
\[
\max_{\ell\in\{k,\dots,n\}}\frac{\la_\ell(\intfloorsm{n^\beta}+1:\ell)}{
\la_\ell}
\ge \frac{c_4}{n^{\beta(2-\alpha)}}
\]
for some $c_4>0$. We suppose that $n$ is large enough
to satisfy $1-n^{-\gamma}\ge 1/2$. These observations,
Lemma~\ref{lem:latail} and Corollary~\ref{cor:lan}
allow us to rewrite inequality~\eqref{eq:bounds:t5} as
\begin{align*}
&\left(\frac{A\Gamma(\alpha+1)}{2-\alpha}n^{2-\alpha}
+O\left(n^{2-\sigma-\varsigma'}\right)\right)
\left(\frac{A\alpha}{2-\alpha}\,\frac{\Gamma(m+\alpha-2)}{\Gamma(m)}
b^{2-\alpha}
+O\left(\frac{\Gamma(m+\alpha+\varsigma'-2)}{\Gamma(m)}
b^{2-\alpha-\varsigma'}\right)\right)
\\
&\qquad\qquad\le 
\left(\frac{A\alpha}{2-\alpha}\,\frac{\Gamma(m+\alpha-2)}{\Gamma(m)}
n^{2-\alpha}
+O\left(\frac{\Gamma(m+\alpha+\varsigma'-2)}{\Gamma(m)}n^{
2-\alpha-\varsigma' }\right)\right)\\
&\qquad\qquad\qquad\qquad\qquad
\times\left(\frac{A\Gamma(\alpha+1)}{2-\alpha}b^{2-\alpha}
+O\left(b^{2-\alpha-\varsigma'}\right)\right)
\left(1+n^{-\gamma}\right)+c_5 b^{2-\alpha}n^{(1-\beta)(2-\alpha)}
\end{align*}
for some $c_5>0$.
Simplification shows that this inequality holds provided
\[
c_6\frac{\Gamma(m+\alpha+\varsigma'-2)}{\Gamma(m)}b^{2-\alpha-\varsigma'
}
n^ { 2-\alpha }
\le 
\frac{\Gamma(m+\alpha-2)}{\Gamma(m)}b^{2-\alpha}n^{2-\alpha-\gamma}
+c_7 b^{2-\alpha}n^{(1-\beta)(2-\alpha)}
\]
for suitable constants $c_6,c_7>0$.  Further simplification gives
\begin{equation}\label{eq:bounds:t7}
c_6\le
b^{\varsigma'}\left(\frac{\Gamma(m+\alpha-2)}{
\Gamma(m+\alpha+\varsigma'-2)}n^{-\gamma}
+c_7\frac{\Gamma(m)}{\Gamma(m+\alpha+\varsigma'-2)}n^{-\beta(2-\alpha)}
\right)
\,.
\end{equation}
Proceeding as above, the expression in brackets attains its minimum in
$m\in{[3,n^\beta]}$
at $m''\sim c_8 n^{\beta-\gamma/(2-\alpha)}$, $c_8>0$, 
with the minimum value asymptotic to
\[
\frac{\Gamma(m''+\alpha-2)}{\Gamma(m''+\alpha+\varsigma'-2)}n^{-\gamma}
+c_7\frac{\Gamma(m'')}{\Gamma(m''+\alpha+\varsigma'-2)}
n^{-\beta(2-\alpha) }
\sim c_9 n^{-\beta\varsigma'-\gamma(2-\alpha-\varsigma')/(2-\alpha)}
\]
where $c_9>0$.
Since $b\ge \intfloor{n^\upsilon}$ the right-hand side
of~\eqref{eq:bounds:t7}
grows to infinity as $n\to\infty$ as long as \eqref{eq:ineqboundsexist}
holds.
This observation finishes the proof.
\end{proof}

We want to keep control over the difference between distributions
of $J^+_{n\fromto k}$, $J^-_{n\fromto k}$ and $J_b$, $n\ge b\ge k$.
In particular, the following statement provides bounds for
divergence of means.

\begin{lemma}\label{lem:expm} Suppose $\beta<\upsilon<1$. Then there
exists $c>0$ such that for $n$ large enough and for $k$ in range $n\ge
k\ge \intfloor{n^\upsilon}$ the following inequalities hold:
\begin{align*}
\left|\ex\bigl[J^-_{n\fromto k}-1\bigr]-\ex\bigl[J_n-1\bigr]\right|
&{}\le c\,\max\left\{n^{-\gamma},
  n^{-\beta(1-\alpha)},
  n^{\theta-\beta(2-\alpha)},
  n^{1-\theta(2-\alpha)}\right\}\,,
\\
\left|\ex\bigl[J^+_{n\fromto k}-1\bigr]-\ex\bigl[J_k-1\bigr]\right|
&{}\le c\,\max\left\{n^{-\beta(1-\alpha)},n^{-\gamma}\right\}\,.
\end{align*}
\end{lemma}

\begin{proof}
We start with the following observation.
For $2\le m< b$, as
$b\to\infty$ but $m/b\to 0$,
 \begin{align*}
 \sum_{j=m}^b j\la_{b,j}
 &{} =\int_0^1\sum_{j=m}^bj\binom{b}{j} x^{j-2}(1-x)^{b-j}\Lambda(dx)
=-\int_0^1\sum_{j=m}^bj\binom{b}{j} x^{j-1}(1-x)^{b-j} dG_{-1}(x)\\
 &{}=m(m-1)\binom{b}{m}\int_0^1 x^{m-2}(1-x)^{b-m}G_{-1}(x)\,dx
  \sim \frac{A\alpha\Gamma(b+1)\Gamma(m+\alpha-2)}
  {(1-\alpha)\Gamma(b+\alpha-1)\Gamma(m-1)}\,,
 \end{align*}
 where 
 \[
 G_{-1}(x)=\int_x^1\frac{\Lambda(dy)}{y} 
 \sim \frac{A\alpha}{1-\alpha}x^{\alpha-1},~~~{\rm as~~~} x\to0\,.
 \]
Taking $m=\intfloor{n^\beta}+1$, for some $\beta\in{]0,\upsilon[}$
we see using Corollary~\ref{cor:lan} that
\begin{equation}\label{eq:meantail}
\sum_{j=\intfloorsm{n^\beta}+1}^k \frac{j\la_{k,j}}{\la_k}
\sim \frac{2-\alpha}{(1-\alpha)\Gamma(\alpha)}\,n^{-\beta(1-\alpha)}
\end{equation}
as $n,k\to\infty$ with $n\ge k\ge \intfloor{n^\upsilon}$.

Now the proof follows by a simple calculation.  The mean of
$J^-_{n\fromto k}$ can be estimated using Lemma~\ref{lem:exJn}
and~\eqref{eq:meantail}:
\begin{equation*}
\begin{split}
\ex\bigl[J^-_{n\fromto k}-1\bigr]
&{}=\frac{\la_{n,2}-n^{-\gamma}\la_n(3:\intfloorsm{n^\beta})
 +\la_n(\intfloorsm{n^\beta}:n)}{\la_n}
+\sum_{j=3}^{\intfloorsm{n^\beta}}
  \frac{(j-1)\la_{n,j}(1+n^{-\gamma})}{\la_n}\\
+& 2\left(\intfloor{n^\theta}-2\right)
\max_{\ell\in\{k,\dots,n\}}
  \frac{\la_\ell(\intfloorsm{ n^\beta } +1:\ell)}{\la_\ell}
+2\left(n-\intfloor{n^\theta}\right)
\max_{\ell\in\{k,\dots,n\}}
  \frac{\la_\ell(\intfloorsm{ n^\theta } +1:\ell)}{\la_\ell}\\
&{}=\ex\bigl[J_n-1\bigr]\left(1+n^{-\gamma}\right)
-\frac{n^{-\gamma}\la_n(2:\intfloorsm{n^\beta})}{\la_n}
+\frac{\la_n(\intfloorsm{n^\beta}+1:n)}{\la_n}\\
&{}\qquad-\sum_{j=\intfloorsm{n^\beta}+1}^n
  \frac{(j-1)\la_{n,j}}{\la_n}
+O\left(\max\left\{n^{\theta-\beta(2-\alpha)},
  n^{1-\theta(2-\alpha)}\right\} \right)\\
&{}=\ex\bigl[J_n-1\bigr]
+O\left(\max\left\{n^{-\gamma},
  n^{-\beta(1-\alpha)},
  n^{\theta-\beta(2-\alpha)},
  n^{1-\theta(2-\alpha)}\right\} \right)\,.
\end{split}
\end{equation*}
Similarly, since $\upsilon>\beta$ formula~\eqref{eq:meantail} is
applicable and implies together with Lemma~\ref{lem:latail} that
\begin{align*}
\ex\bigl[J_{n\fromto k}^+-1\bigr]
={}&\frac{\la_{k,2}+n^{-\gamma}\la_k(3:\intfloorsm{n^\beta})
+\la_k(\intfloorsm{n^\beta }+1:k)}{\la_k}
+\sum_{j=3}^{\smash[t]{\intfloor{n^\beta}}}
\frac{(j-1)\la_{k,j}}{\la_k}
\left(1-n^{-\gamma}\right)\\
=\ex\bigl[J_k-1\bigr]&{}\left(1-n^{-\gamma}\right)
-\sum_{j=\intfloor{n^\beta}+1}^{k}\frac{(j-1)\la_{k,j}}{\la_k}
  \left(1-n^{-\gamma}\right)
+n^{-\gamma}\frac{\la_k(2:\intfloorsm{n^\beta})}{\la_k}
+\frac{\la_k(\intfloorsm{n^\beta}+1:k)}{\la_k}\\
={}&\ex\bigl[ J_k-1 \bigr]
+O\left(\max\left\{n^{-\beta(1-\alpha)},n^{-\gamma}
\right\}\right)\,,
\end{align*}
so the claim follows.
\end{proof}

Using a familiar device,
Lemma~\ref{lem:bounds} enables us to couple random variables
$J^+_{n\fromto k}$, $J_b$ and $J^-_{n\fromto k}$ 
in such a way that
\begin{equation}\label{eq:JJJJcoupling}
J^+_{n\fromto k}\le J_b \le J^-_{n\fromto k}
\end{equation}
holds almost surely.  From this  Lemma~\ref{lem:Green} will follow by comparing ${\cal M}_n$ with two random walks.

\begin{proof}[Proof of Lemma~\ref{lem:Green}]
For $n\ge b\ge k\ge1$ let $g^+(n\fromto k,b)$ and $g^-(n\fromto k,b)$ be
the Green kernels of decreasing random walks started at $n$ 
with decrements $(J_{n\fromto k}^+-1)$ 
and $(J_{n\fromto k}^--1)$, correspondingly.  Suppose
that~\eqref{eq:JJJJcoupling} holds for some $n\ge k$ and for all $b$
in range $n\ge b\ge k$.
Then we have
$g^+(n\fromto k,b)\ge g(n,b)\ge g^-(n\fromto k,b)$ for all such~$b$. 

Take $\gamma>0$ and
$\upsilon>\theta>\beta>1/(2-\alpha)$.  Combination of
Lemmas~\ref{lem:exJn} and~\ref{lem:expm} 
implies that $\ex\bigl[J^+_{n\fromto k}\bigr]\to 1/(1-\alpha)$ and 
$\ex\bigl[J^-_{n\fromto k}\bigr]\to 1/(1-\alpha)$ as $n,k\to\infty$ with $k$ in
range $n\ge k\ge \intfloor{n^\upsilon}$.

Assume further that parameters $\beta$, $\gamma$, $\upsilon$ and
$\theta$ are chosen to satisfy 
inequality~\eqref{eq:ineqboundsexist}.
Then the coupling~\eqref{eq:JJJJcoupling} exists for $n\ge b\ge
k\ge\intfloor{n^\upsilon}$ by~Lemma~\ref{lem:bounds}. 
Applying a standard result of renewal theory,
\[
g^+(n\fromto k,b)\sim\frac{1}{\ex\bigl[J^+_{n\fromto k}-1\bigr]}\to 1-\alpha
\qquad\text{and}\qquad
g^-(n\fromto k,b)\sim\frac{1}{\ex\bigl[J^-_{n\fromto k}-1\bigr]}\to 1-\alpha
\]
as $n,k\to\infty$ and $n-b\to\infty$.
Thus $g(n,b)$ has the same limit.
\end{proof}

\section{The total number of collisions}\label{sec:limitthm}

 We are
in position now to present our main result on the convergence of the number of collisions $C_n$
in the $\Lambda$-coalescent on $n$ particles.

\begin{thm}\label{thm:numcoll}
Suppose that the measure $\Lambda$ satisfies~\eqref{ass2} 
with
$\varsigma>    \max\left\{\frac{(2-\alpha)^2}{
5-5\alpha+\alpha^2},1-\alpha\right
\}$. 
Then, as $n\to\infty$, we have the convergence in distribution 
\[
\frac{C_n-(1-\alpha)n}{(1-\alpha)n^{1/(2-\alpha)}}\to_d \CS_{2-\alpha}
\]
 to a stable random variable $\CS_{2-\alpha}$  with the characteristic
function
\begin{equation}\label{eq:charfS2a}
\ex\left[e^{\imi u\CS_{2-\alpha}}\right] =
\exp\left(-e^{-\imi \pi\alpha\sign(u)/2}|u|^{2-\alpha}\right).
\end{equation}
\end{thm}

\noindent 
We emphasize that $\varsigma=1$ satisfies assumptions
of the above Theorem for all $\alpha\in{]0,1[}$. This is important
because $\varsigma=1$ appears, say, if $\Lambda$ is a beta-measure.

\smallskip
\noindent
{\bf Remark.}
The characteristic function~\eqref{eq:charfS2a} is not
a canonic form for the characteristic function of stable distribution.  
There are several commonly used 
parametrisations for stable variables, see 
\cite{Samorodnitsky-Taqqu,Zolotarev}; the difference between them 
being  a frequent source of  confusion.  
Apparently 
the most common parametrisation
involves the
 \textit{index of stability} $\alpha\in{]0,2]}$,
the
\textit{skewness}
 $\beta\in{[-1,1]}$,
the \textit{scale} 
$\sigma>0$ 
and the \textit{location}
 $\mu\in\BR$,
 so that a
random variable
$\CS$
has the stable distribution with parameters $(\alpha,\beta,\sigma,\mu)$
iff
\[
\ex\left[e^{\imi u\CS}\right]=\begin{cases}
\exp\left(-\sigma^\alpha|u|^\alpha\left(1-\imi \beta\tan\frac{\pi\alpha}{2}
\sign u\right)+\imi \mu u\right), 
  \qquad &\alpha\ne1,\\
\exp\left(-\sigma|u| \left(1+\frac{2\imi \beta}{\pi}(\sign
u)\log|u|\right)+\imi \mu u\right), &\alpha=1. 
\end{cases}
\]
In this parametrisation our random variable $\CS_{2-\alpha}$ has 
$\left(2-\alpha,-1,(\cos\frac{\pi\alpha}{2})^{1/(2-\alpha)},0\right)$-stable 
distribution since its charactersitic function can be rewritten as
\begin{align*}
\exp\left(-e^{-\imi \pi\alpha\sign(u)/2}|u|^{2-\alpha}\right)
&{}=\exp\left(-\cos\tfrac{\pi\alpha}{2}|u|^{2-\alpha}(1-\imi \tan\tfrac{
\pi\alpha}{2}\sign u)\right)\\
&{}=\exp\left(-\cos\tfrac{\pi\alpha}{2}|u|^{2-\alpha}(1+\imi \tan\tfrac{
\pi(2-\alpha)}{2}\sign u)\right)\,.
\end{align*}
Thus $\CS_{2-\alpha}$ has $(2-\alpha)$-stable distribution totally
skewed to the left.
\vskip0.5cm

The main idea of the proof is that the decrements~$J_b$  of
are almost identically distributed  
for large $b$, as Corollary~\ref{cor:lan} suggests.  
However the
nonstationarity prevents the possibility of
any direct analysis.  
To override this, we use the technique of
stochastic bounds  introduced in the previous section.  First
we introduce some auxiliary notations.

For $1\le k\le n$ the coalescent started with $n$ particles after some
series of collisions will reach a state with less than $k+1$ particles; 
let $C_{n\fromto k}$ denote
the number of collisions and let
$B_{n,k}\leq k$ denote the number of particles 
as the coalescent enters such state.
In particular,  $C_n=C_{n\fromto1}$. 
For $J_{n\fromto k,m}^\pm$  independent copies of
 $J_{n\fromto k}^\pm$, introduce
\begin{equation}\label{eq:defcpm}
C_{n\fromto k,\ell}^+ :=\min\left\{c:
\sum_{m=1}^c \left(J_{n\fromto k,m}^+-1\right)\ge \ell\right\}\,,
\qquad  
C_{n\fromto k,\ell}^-:=\min\left\{c:
\sum_{m=1}^c \left(J_{n\fromto k,m}^--1\right)\ge \ell\right\},
\end{equation}
the minimal number of decrements distributed as 
$J_{n\fromto k}^+-1$ (respectively, 
$J_{n\fromto k}^--1$) needed
to drop by at least $\ell$.   
We skip the 
index $\ell$ when it is equal to $n-k$, so 
that $C^\pm_{n\fromto k}\equiv C^\pm_{n\fromto k,n-k}\,$.

Under assumptions of Lemma~\ref{lem:bounds} we can couple the 
corresponding Markov chains so that \eqref{eq:JJJJcoupling} holds
almost surely for all large enough $n$ once $n\ge b\ge k\ge
\intfloor{n^\upsilon}$. Consequently, for such $n$  the coupled
Markov chains satisfy
\[
C_{n\fromto \intfloor{n^{\upsilon}}}^+\ge C_{n\fromto
\intfloor{n^{\upsilon}}}\ge
C_{n\fromto\intfloor{n^{\upsilon}}}^-\,.
\]  
In other words,
\begin{equation}\label{ineqCCC}
C_{n\fromto \intfloor{n^{\upsilon}}}^+\ge_d
C_{n\fromto \intfloor{n^{\upsilon}}}\ge_d
C_{n\fromto \intfloor{n^{\upsilon}}}^-\,.
\end{equation}

In order to find the limit distributions for
$C_{n\fromto\intfloor{n^\upsilon}}^\pm$ we
need the following statement about the characteristic function 
\[
\phi_n(u):=\ex\left[e^{\imi u(J_n-1)}\right]
\]
of the first decrement  of ${\cal M}_n$. 

\begin{lemma}\label{lem:charJ} 
Let $\Lambda$ satisfy \eqref{ass2} 
with $\varsigma> 1-\alpha$. 
Then there exists $\delta>0$ such that  
\[
\phi_n\left(s/m\right)
=1+\frac{\imi s}{(1-\alpha)m}
-\frac{\omega(s)|s|^{2-\alpha}}{(1-\alpha)m^{2-\alpha}}
+O\left(m^{\alpha-2-\delta}\right)\,,
\]
as $n,m\to\infty$ with $m\le n^{\upsilon}$ for some $\upsilon<1$, 
where $\omega(s)=e^{\imi \pi\alpha\sign(s)/2}$.  
\end{lemma}

\begin{proof}
We write for shorthand $u=s/m$.
For $u=0$ the claim is obvious, so we suppose that $u\ne0$.
The characteristic function of $J_n-1$ can be written in terms of
$\Lambda$ as
follows:
\[
\phi_n(u)=e^{-\imi u}\sum_{j=2}^n\frac{\la_{n,j}}{\la_n}e^{\imi ju}
=\frac{e^{-\imi u}}{\la_n}\int_0^1
\frac{\left(1-(1-e^{\imi u})x\right)^n-(1-x)^n-nxe^{\imi u}(1-x)^{n-1}}{x^2}
\Lambda(dx)
\]
using the integral representation of $\la_{n,j}$. Denote the
numerator of
the fraction under the integral above by $h_n(u,x)$; then
\[
h_{n+1}(u,x)-h_n(u,x)=x(1-e^{\imi u})\left((1-x)^n-(1-(1-e^{\imi u})x)^n\right)
+x^2n(1-x)^{n-1}e^{\imi u}
\]
so using \eqref{lan} we obtain  
\begin{equation}\label{IntePhi}
\phi_n(u)=1-\frac{1-e^{-\imi u}}{\la_n}\sum_{j=1}^{n-1}\int_0^1
  \frac{(1-x)^j-(1-(1-e^{\imi u})x)^j}{x}\Lambda(dx)
\end{equation}
because $h_1(u,x)=0$.
Taking again differences of $(1-x)^j-(1-(1-e^{\imi u})x)^j$ with respect to
$j$
and calculating it directly for $j=0$ we represent the integral in (\ref{IntePhi})  as
\[
(1-e^{\imi u})\sum_{k=0}^{j-1}\int_0^1(1-(1-e^{\imi u})x)^k\Lambda(dx)
  -\sum_{k=0}^{j-1}\int_0^1(1-x)^k\Lambda(dx).
\]
Exchanging the sums and utilising notation~\eqref{eq:defnu}
for moments $\nu_k$ of~$\Lambda$ we get 
\begin{equation}\label{eq:charJ:1}
\phi_n(u)=1
  +\frac{(1-e^{-\imi u})}{\la_n}\sum_{k=0}^{n-2}(n-k-1)\nu_k
  +\frac{(1-e^{\imi u})^2e^{-\imi u}}{\la_n}\sum_{k=0}^{n-2}(n-k-1)
  \int_0^1\left(1-(1-e^{\imi u})x\right)^k\Lambda(dx).
\end{equation}
By~\eqref{eq:expjn} and 
Lemma~\ref{lem:exJn} the second term above is
\[
\frac{(1-e^{-\imi u})}{\la_n}\sum_{k=0}^{n-2}(n-k-1)\nu_k
=(1-e^{-\imi u})\ex\bigl[J_n-1\bigr]=\frac{1-e^{-\imi u}}{1-\alpha}
  \left(1+O\left(n^{\alpha-1}\right)\right)
\]
since $\varsigma>1-\alpha$ by hypothesis.
Recalling notation $u=s/m$ and inequality $n\ge m^{1/\upsilon}$ 
with $\upsilon<1$ we see that 
\[
\frac{(1-e^{-\imi s/m})}{\la_n}\sum_{k=0}^{n-2}(n-k-1)\nu_k =
\frac{\imi s}{(1-\alpha)m}+O\left(m^{\alpha-2-\delta_1}\right)
\]
for some $\delta_1>0$.
Thus it remains to estimate the last summand in~\eqref{eq:charJ:1}.

Integration by parts  gives
\[
\int_0^1\left(1-(1-e^{\imi u})x\right)^k\Lambda(dx)=
e^{\imi ku}+k(1-e^{\imi u})\int_0^1(1-(1-e^{\imi u})x)^{k-1}\Lambda[0,x] dx.
\]
Substitution of this relation into~\eqref{eq:charJ:1} leads to  
\begin{align}\notag
&\phi_n\left(s/m\right)=1+\frac{\imi s}{(1-\alpha)m}
+\frac{(1-e^{\imi s/m})^2e^{-\imi s/m}}{\la_n}\sum_{k=0}^{n-2}(n-k-1)e^{\imi ks/m}
\\ 
\label{eq:charJ:2}
&{}+\frac{(1-e^{\imi s/m})^3e^{-\imi s/m}}{\la_n}\int_{0}^1\sum_{k=0}^{n-2}
k(n-k-1)
  (1-(1-e^{\imi s/m})x)^{k-1}\Lambda[0,x] dx
+O(m^{-1-\varsigma/\upsilon})\,.
\end{align}
Summation yields
\[
\frac{(1-e^{\imi s/m})^2e^{-\imi s/m}}{\la_n}\sum_{k=0}^{n-2}(n-k-1)e^{\imi ks/m} 
=\frac{e^{-\imi s/m}(n(1-e^{\imi s/m})+e^{\imi sn/m}-1)}{\la_n}\,.
\]
For $m$ big enough and $n\ge m^{1/\upsilon}$ with $\upsilon<1$
\[
\left|\frac{n(1-e^{\imi s/m}) + e^{\imi ns/m}-1}{\la_n}\right|\le
\frac{n|1-e^{\imi s/m}|+|e^{\imi ns/m}-1|}{\la_n} \le
\frac{\text{const}}{n^{1-\alpha}m}=O\left(m^{\alpha-2-\delta_2}\right)
\]
for some $\delta_2>0$.

Let $\theta\in[-\pi/2,\pi/2]$ be such that
$e^{\imi \theta}=\frac{1-e^{\imi s/m}}{|1-e^{\imi s/m}|}$. Note that
$\theta=-\pi\sign(s)/2+O(1/m)$ as $m\to\infty$.
For any $\beta>0$ we have
\begin{multline*} 
\int_0^{1}
  k(1-e^{\imi s/m})(1-(1-e^{\imi s/m})x)^{k-1}x^\beta dx\\
=e^{\imi \theta}\int_0^{k|1-e^{\imi s/m}|}\left(1-e^{\imi \theta}t/k\right)^{k-1} 
  \frac{t^\beta}{(k|1-e^{\imi s/m}|)^\beta}\,dt
= \frac{e^{-\imi \beta\theta}}{\left(k|s|\big/m\right)^\beta} \Gamma(\beta+1)
 \left(1+O(1/m)\right)
\end{multline*}
as $m,k\to\infty$ with $k\ge m^{1+\delta_3}$ for any $\delta_3>0$. 
By assumption~\eqref{ass2} we can write 
$\Lambda[0,x]=Ax^\alpha+f(x)$ where $|f(x)|\le c x^{\alpha+\varsigma}$
for some $c>0$ and all $x\in[0,1]$. Thus, as $m,k\to\infty$ with 
$k\ge m^{1+\delta_3}$,
\[
\int_0^{1} k(1-e^{\imi s/m})(1-(1-e^{\imi s/m})x)^{k-1}\Lambda[0,x]\,dx
 =\frac{Ae^{\imi \pi\alpha\sign(s)/2}}{\left(k|s|\big/ m\right)^\alpha}
\Gamma(\alpha+1)
 +O\left(\frac{m^{\alpha+\varsigma}}
 {k^{\alpha+\varsigma}}+\frac{1}{m^{1-\alpha }k^\alpha}\right).
\]
Take $\delta_3=(1/\upsilon-1)/2$ and denote
$n_0=\intfloor{m^{1+\delta_3}}$.
Divide the last sum 
in~\eqref{eq:charJ:2} into two sums over 
$k\ge n_0$ and $k<n_0$.  The first sum
is estimated taking \eqref{obiglan} into 
account as
\begin{align*}
\frac{(1-e^{\imi s/m})^2}{\la_n}&\sum_{k=n_0}^{n-2}(n-k-1)
  \int_0^{1}  k(1-e^{\imi s/m})(1-(1-e^{\imi s/m})x)^{k-1}\Lambda[0,x]\,dx\\
&= -\frac{|s|^{2-\alpha}e^{\imi \pi\alpha\sign(s)/2}(2-\alpha)}
  {m^{2-\alpha}n^{2-\alpha}}
  \sum_{k=n_0}^{n-2} (n-k)k^{-\alpha} \\
&\qquad +O\left(\frac{1}{n^{2-\alpha}m^{2-\alpha-\varsigma}}
  \sum_{k=n_0}^{n-2}\frac{n-k}{k^{\alpha+\varsigma}} 
+\frac{1}{n^{2-\alpha}m^{3-\alpha}}\sum_{k=n_0}^{n-2}\frac{n-k}{k^\alpha
}   \right) \\
&=
-\frac{|s|^{2-\alpha}e^{\imi \pi\alpha\sign(s)/2}(2-\alpha)}{m^{2-\alpha}}\,
  \int_{m^{1+\delta_3}n^{-1}}^1 
x^{-\alpha}(1-x)dx+O\left(\frac{1}{n^{1-\alpha}m^{1+\varsigma\delta_3}}
+\frac{1}{m^{3-\alpha}}\right)\\
& 
=-\frac{|s|^{2-\alpha}e^{\imi \pi\alpha\sign(s)/2}}{(1-\alpha)m^{2-\alpha}}
  +O(m^{\alpha-2-\delta_4})
\end{align*}
for some $\delta_4>0$.  The same argument
applied to the sum over  $k=0,\dots,n_0-1$ shows
that it constitutes a lower order term
to the whole sum.  Thus it remains to combine the results above to get
the statement of Lemma.  
\end{proof}

Next we show that under certain
assumptions the same asymptotic expansion is also valid for the
characteristic functions of $J^\pm_{n\fromto k}$
\[
\phi^+_{n\fromto k}(u):=\ex\left[e^{\imi uJ^+_{n\fromto k}}\right]
\qquad\text{ and }\qquad
\phi^-_{n\fromto k}(u):=\ex\left[e^{\imi uJ^-_{n\fromto k}}\right]\,.
\]

\begin{lemma}\label{lem:charJpm}
Suppose~\eqref{ass2} holds with
$\varsigma>1-\alpha$ and that 
the parameters in~\eqref{eq:defqpm} satisfy
inequalities
\begin{equation}\label{eq:ineqcharJpm}
\gamma>\frac{1-\alpha}{2-\alpha},\qquad \text{and}\qquad
\upsilon>\theta>\beta>\frac{1}{2-\alpha}.
\end{equation}
Then there exists $\delta>0$ such that
\begin{align*}
\phi^+_{n\fromto k}\left(s/m\right)
&{}=1+\frac{\imi s}{(1-\alpha)m}
-\frac{\omega(s)|s|^{2-\alpha}}{(1-\alpha)m^{2-\alpha}}
+O\left(m^{\alpha-2-\delta}\right)\,,\\
\phi^-_{n\fromto k}\left(s/m\right)
&{}=1+\frac{\imi s}{(1-\alpha)m}
-\frac{\omega(s)|s|^{2-\alpha}}{(1-\alpha)m^{2-\alpha}}
+O\left(m^{\alpha-2-\delta}\right)\,,
\end{align*}
as $n,k,m\to\infty$ in such a way that $n\ge k\ge \intfloor{n^\upsilon}$
and $m\le c n^{1/(2-\alpha)}$ for some $c>0$.
\end{lemma}

\begin{proof} 
The characteristic function of $J_{n\fromto k}^--1$ is by definition
\begin{multline*}
\phi_{n\fromto k}^-(u)
=e^{\imi u}\left(\frac{\la_{n,2}-n^{-\gamma}\la_n(3:\intfloorsm{
n^\beta})
+\la_n(\intfloorsm{n^\beta}+1:n)}{\la_n}-2\max_{\ell\in\{k,\dots,n\}}
\frac { \la_\ell(\intfloorsm { n^\beta } +1:\ell)}{\la_\ell}\right)\\
+\sum_{j=3}^{\intfloor{n^\beta}}\frac{\la_{n,j}(1+n^{-\gamma})}{
\la _n} e^{\imi (j-1)u} 
+ 2e^{\imi (\intfloorsm{n^\theta}-1)u}
\max_{\ell\in\{k,\dots,n\}}
  \frac{\la_\ell(\intfloorsm{ n^\beta } +1:\ell)}{\la_\ell}\\
+ 2\left(e^{\imi (n-1)u}-e^{\imi (\intfloorsm{n^\theta}-1)u}\right)
\max_{\ell\in\{k,\dots,n\}}
  \frac{\la_\ell(\intfloorsm{ n^\theta } +1:\ell)}{\la_\ell}\,.
\end{multline*}
We rewrite it as follows:
\begin{multline*}
\phi_{n\fromto k}^-(u)=\phi_n(u)\left(1+n^{-\gamma}\right)
-n^{-\gamma}e^{\imi u}
+\frac{\la_n(\intfloorsm{n^\beta}+1:n)(1+n^{-\gamma}) } {\la_n} e^{\imi u}\\
-\sum_{j=\intfloorsm{n^\beta}+1}^n\frac{(1+n^{-\gamma})\la_{n,j}}
{ \la_n} e^ { \imi (j-1)u }
+ 2\left(e^{\imi (\intfloorsm{n^\theta}-1)u}-e^{\imi u}\right)
\max_{\ell\in\{k,\dots,n\}}
  \frac{\la_\ell(\intfloorsm{ n^\beta } +1:\ell)}{\la_\ell}\\
+ 2\left(e^{\imi (n-1)u}-e^{\imi (\intfloorsm{n^\theta}-1)u}\right)
\max_{\ell\in\{k,\dots,n\}}
  \frac{\la_\ell(\intfloorsm{ n^\theta } +1:\ell)}{\la_\ell}\,.
\end{multline*}
Four last summands are of the order of $n^{-\beta(2-\alpha)}$ by
Lemma~\ref{lem:latail}.  From~\eqref{eq:ineqcharJpm},
$\beta>1/(2-\alpha)$ and so the bound $m\le c n^{1/(2-\alpha)}$
guarantees that these four summands constitute
$O(m^{\alpha-2-\delta_1})$ terms to the whole sum, for
$\delta_1,\delta_2,\dots$ some positive constants. 
The same bound on $m$ allows application of
Lemma~\ref{lem:charJ} for $\phi_n(s/m)$ which leads to
\begin{align*}
\phi^-_{n\fromto k}(s/m)
&{}=\left(1+\frac{\imi s}{(1-\alpha)m}-\frac{\omega(s)|s|^{2-\alpha}}{
(1-\alpha)m^{2-\alpha}}+O(m^{\alpha-2-\delta_2})\right)\left(1+n^{
-\gamma } \right)\\
&\qquad-n^{-\gamma}\left(1+\frac{\imi s}{m}+O(m^{-2})\right)
+O(m^{\alpha-2-\delta_1})\\ 
&{}=1+\frac{\imi s}{(1-\alpha)m}-\frac{\omega(s)|s|^{2-\alpha}}{
(1-\alpha)m^{2-\alpha}}+\left(\frac{1}{1-\alpha}-1\right)
  \frac{\imi s\, n^{-\gamma}}{m}
+O(m^{\alpha-2-\delta_3}) 
\end{align*}
and the claim about $\phi^-_{n\fromto k}$ follows from inequality
$\gamma>(1-\alpha)/(2-\alpha)$.

Analogously,
\begin{align*}
\phi_{n\fromto k}^+(u)
&{}=\frac{\la_{k,2}+n^{-\gamma}\la_k(3:\intfloorsm{n^\beta})
+\la_k(\intfloorsm{n^\beta}+1:k)} { \la_k } e^ { \imi u }
+\sum_{j=3}^{\intfloorsm{n^\beta}}
\frac{\la_{k,j}(1-n^{-\gamma})e^{\imi (j-1)u}}{\la_k}\\
&{}=\frac{\la_{k,2}}{\la_k}\left(1-n^{-\gamma}\right)e^{\imi u}
+n^{-\gamma}e^{\imi u}+\frac{\la_k(\intfloorsm{n^\beta}+1:b)}{\la_k}
\left(1-n^{-\gamma}\right)e^{\imi u}\\
&\qquad\qquad\qquad+\sum_{j=3}^n\frac{\la_{k,j}(1-n^{-\gamma})}
{\la_k}e^{ \imi (j-1)u }
-\sum_{j=\intfloorsm{n^\beta}+1}^n\frac{\la_{k,j}(1-n^{-\gamma})}
  { \la_k } e^ {\imi (j-1)u}\\
&{}=\phi_k(u)
\left(1-n^{-\gamma}\right)+n^{-\gamma}e^{\imi u}
+\left(1-n^{-\gamma}\right)\sum_{j=\intfloorsm{n^\beta}+1}^n
\frac{\la_{k,j}(1-e^{\imi ju})e^{\imi u} }
  { \la_k } \,.
\end{align*}
Since $k\ge \intfloor{n^\upsilon}$ with $\upsilon>\beta>1/(2-\alpha)$
and $m\le cn^{1/(2-\alpha)}$, Lemma~\ref{lem:latail} ensures that the
last sum above is $O(n^{-\beta(2-\alpha)})=O(m^{\alpha-2-\delta_4})$. 
Since $m$ grows slower than $k$ by hypothesis, Lemma~\ref{lem:charJ}
can be applied for $\phi_k(s/m)$ and the claim again follows from
inequality $\gamma>(1-\alpha)/(2-\alpha)$.
\end{proof}

\begin{lemma}\label{lem:limit4sum}
Suppose~\eqref{ass2} holds with
$\varsigma>1-\alpha$ and that 
the parameters in~\eqref{eq:defqpm} satisfy
inequalities
\begin{equation}\label{eq:ineqsamesum}
\gamma>\frac{1-\alpha}{2-\alpha},\qquad 
1>\upsilon>\theta>\beta>\frac{5-5\alpha+\alpha^2}{(2-\alpha)^3},\qquad
\beta(2-\alpha)-\frac{1-\alpha}{2-\alpha}>\theta>\frac{3-2\alpha}{
(2-\alpha)^2 }\,.
\end{equation}
Let $S_{n\fromto k,h}^+$, respectively $S_{n\fromto
k,h}^-$,  be the sum of\/ $h$ independent copies of\/
$J_{n\fromto k}^+-1$, respectively  $J_{n\fromto k}^--1$.  
Then 
\[
\frac{S_{n\fromto k,h}^+-h/(1-\alpha)}
{\left(h/(1-\alpha)\right)^{1/(2-\alpha)}}
\to_d \bar\CS_{2-\alpha}
\qquad\text{and}\qquad
\frac{S_{n\fromto k,h}^--h/(1-\alpha)}
{\left(h/(1-\alpha)\right)^{1/(2-\alpha)}}
\to_d \bar\CS_{2-\alpha}
\]
as $n,h\to\infty$ with $h=O(n)$ and $n\ge k\ge
\intfloor{n^\upsilon}$, $1>\upsilon>\beta$, where
$\bar\CS_{2-\alpha}$ is a stable
random variable with the characteristic function
\begin{equation}\label{eq:charfbarS2a}
\ex\left[e^{\imi u\bar\CS_{2-\alpha}}\right]=\exp
\left[-\omega(u)|u|^{2-\alpha}\right]\,,\qquad\qquad\omega(u)=\exp\left(\frac{
\imi \pi\alpha\sign u}{2}\right).
\end{equation}
\end{lemma}

\begin{proof} First note that a solution of 
inequality~\eqref{eq:ineqsamesum}
always exists. Since~\eqref{eq:ineqcharJpm} follows
from~\eqref{eq:ineqsamesum}, the bound $h=O(n)$ guarantees that
Lemma~\ref{lem:charJpm} is applicable with
$m=\left(h/(1-\alpha)\right)^{1/(2-\alpha)}$.
Lemmas~\ref{lem:exJn} and~\ref{lem:expm} provide tough bounds for
$\ex\bigl[J^\pm_{n\fromto k}-1\bigr]$. Namely, inequality~\eqref{eq:ineqsamesum}
implies that
\begin{equation}\label{eq:l4s:t1}
\left|\ex\bigl[J^-_{n\fromto k}-1\bigr]-\frac{1}{1-\alpha}\right|
=O(n^{-(1-\alpha)/(2-\alpha)-\delta_1})
\end{equation}
for some $\delta_1>0$. Hence
\[
\phi_{n\fromto k}^-
\left(\frac{s}{\left(h/(1-\alpha)\right)^{1/(2-\alpha)}}\right) 
\exp\left(-\frac{\ex\bigl[J_{n\fromto k}^--1\bigr]}
{\left(h/(1-\alpha)\right)^{1/(2-\alpha)}}\imi s\right)
=1-\frac{\omega(s)|s|^{2-\alpha}}{h}+O(h^{-1-\delta_2})
\]
for some $\delta_2>0$.
Moreover, equation~\eqref{eq:l4s:t1} and $h=O(n)$ imply
\[
\frac{\ex\bigl[S_{n\fromto k,h}^-\bigr]
-\tfrac{h}{1-\alpha}}{\left(h/(1-\alpha)\right)^{1/(2-\alpha)}}
=O\left(h^{1-1/(2-\alpha)}n^{-(1-\alpha)/(2-\alpha)-\delta_1}
\right)
=O\left(h^{-\delta_1}\right)\,,
\]
as $n,h\to\infty$. Hence, for some $\delta_3>0$ 
\begin{align*}
\ex&\left[\exp
\left(\frac{S_{n\fromto
k,h}^--\frac{h}{1-\alpha}}{\left(h/(1-\alpha)\right)^{1/(2-\alpha)}}
\imi s\right)\right]\\
&\qquad\qquad=\exp\left(
\frac{\ex\bigl[S_{n\fromto k,h}^-\bigr]
-\frac{h}{1-\alpha}}{\left(h/(1-\alpha)\right)^{1/(2-\alpha)}}\imi s
\right)
\ex\left[\exp
\left(\frac{S_{n\fromto k,h}^--\ex\bigl[S_{n\fromto k,h}^-\bigr]}
{\left(h/(1-\alpha)\right)^{1/(2-\alpha)}}\imi s\right)\right]\\
&\qquad\qquad=\left(1+O(h^{-\delta_3})\right)\left(1-\frac{\omega(s)|s|^{2-\alpha}}{
h}+O(h^{-1-\delta_2})\right)\Big.^{\!h}\\
&\qquad\qquad\to \exp\left[ -e^{\frac{\imi \pi\alpha\sign s}{2}}|s|^{2-\alpha}\right]
\end{align*}
and the claim about $S^-_{n\fromto k,h}$ follows.  

Treatment of the limit theorem for $S^+_{n\fromto k,h}$ literally
repeats the above steps and is omitted.
\end{proof}

Let $F_{2-\alpha}(\cdot)$ be the  distribution function
of the stable random variable $\CS_{2-\alpha}$ defined
by~\eqref{eq:charfS2a} 
and $\bar F_{2-\alpha}(\cdot)$ be that of $\bar\CS_{2-\alpha}$ defined 
by~\eqref{eq:charfbarS2a}. 
Note that
the random variables $\CS_{2-\alpha}$ and
$-\bar\CS_{2-\alpha}$ 
have the same distributions, i.e.\ $F_{2-\alpha}(t)=1-\bar F_{2-\alpha}(-t)$.

\begin{lemma}\label{lem:Cntondelta} Let the measure $\Lambda$
satisfy~\eqref{ass2} with
$\varsigma>\max\left\{\frac{(2-\alpha)^2}{
5-5\alpha+\alpha^2},1-\alpha\right\}$.  Then there exists $\upsilon<1$
such that
\[
\frac{C_{n\fromto\intfloor{n^\upsilon}}-(n-\intfloor{
n^\upsilon})(1-\alpha)}
{n^{1/(2-\alpha)}(1-\alpha)}
\to_d {\cal S}_{2-\alpha}\,~~~{\rm as}~~~ n\to\infty.
\]
\end{lemma}

\begin{proof}
Suppose that $\beta$, $\gamma$ and $\theta$ satisfy both inequalities
\eqref{eq:ineqboundsexist} and~\eqref{eq:ineqsamesum}.
This is always possible if $\varsigma$ satisfy condition stated in the
Lemma.  Indeed, the only constraints which can become inconsistent by
joining inequalities~\eqref{eq:ineqboundsexist}
and~\eqref{eq:ineqsamesum} are the constraints on $\gamma$ 
\begin{equation}\label{eq:ineqgamma}
\frac{(\upsilon-\beta)(2-\alpha)\varsigma'}{2-\alpha-\varsigma'}
>\gamma>\frac{1-\alpha}{2-\alpha}\,.
\end{equation}
By~\eqref{eq:ineqsamesum} we can choose $\beta$ and $\upsilon$ such that
$\upsilon-\beta<1-\frac{5-5\alpha+\alpha^2}{(2-\alpha)^3}\,$. 
Hence~\eqref{eq:ineqgamma} is solvable for 
\[
\left(1-\frac{5-5\alpha+\alpha^2}{(2-\alpha)^3}\right)
\frac{(2-\alpha)\varsigma' }{2-\alpha-\varsigma'}
>\frac{1-\alpha}{2-\alpha}\,.
\]
Resolving $\varsigma'$ from this inequality and recalling its
definition leads to the lower bound on
$\varsigma$ in the claim.

By definition $C_{n\fromto\intfloor{n^\upsilon},\,d}^+$ is the random 
number of decrements $J_{n\fromto \intfloor{n^\upsilon}}^+-1$ needed
to make a total move larger than $d$.  Hence for all
$h>0$
\begin{equation}\label{eq:linkCS}
\prob\left[C_{n\fromto \intfloor{n^\upsilon},\,d}^+\le h\right] =
\prob\left[S_{n\fromto \intfloor{n^\upsilon},h}^+\ge d\,\right].
\end{equation}
Take now $d=n-\intfloor{n^\upsilon}$ and
$h=\intfloor{\left(n-\intfloor{n^\upsilon}+t_nn^{1/(2-\alpha)}
\right)(1-\alpha)}$ where $t_n\to t$ as $n\to\infty$.  
Then $h\to\infty$ but $h=O(n)$ as $n\to\infty$. Moreover, 
\[
d=n-\intfloor{n^\upsilon}=\frac{h}{1-\alpha}
-t\left(\frac{h}{1-\alpha}\right)^{1/(2-\alpha)}(1+o(1)),
\qquad n\to\infty,
\] 
and application of Lemma~\ref{lem:limit4sum} ensures that
the RHS of~\eqref{eq:linkCS} converges to $\bar
F_{2-\alpha}(-t)$ as
$n\to\infty$, since $\bar F_{2-\alpha}$ is continuous~\cite{Zolotarev}. 
Thus 
\begin{multline*}
\prob\left[\frac{C_{n\fromto\intfloor{n^\upsilon}}^+-(n-\intfloor{
n^\upsilon})(1-\alpha)}
{n^{1/(2-\alpha)}(1-\alpha)}\le t\right]
\sim\prob\left[\,C_{n\fromto \intfloor{n^\upsilon}}^+\le h\,\right]\\
=\prob\left[\,S_{n\fromto \intfloor{n^\upsilon},h}^+\ge
\frac{h}{1-\alpha}-t\left(\frac{h}{1-\alpha}\right)^{1/(2-\alpha)}
\left(1+o(1)\right)\right]
\to 1-\bar F_{2-\alpha}(-t)=F_{2-\alpha}(t)\,.
\end{multline*}

Replacing $S_{n\fromto \intfloor{n^\upsilon},h}^+$ with 
$S_{n\fromto \intfloor{n^\upsilon},h}^-$ and
$C_{n\fromto\intfloor{n^\upsilon}}^+$ with
$C_{n\fromto\intfloor{n^\upsilon}}^-$ in the above argument we obtain 
\[
\prob\left[\frac{C_{n\fromto\intfloor{n^\upsilon}}^--(n-\intfloor{
n^\upsilon})(1-\alpha)}
{n^{1/(2-\alpha)}(1-\alpha)}\le
t\right]\to  F_{2-\alpha}(t).
\]
Hence the claim follows from inequalities~\eqref{ineqCCC}
since $F_{2-\alpha}$ is continuous.
\end{proof}

\begin{proof}[Proof of Theorem~\ref{thm:numcoll}]
Recall that $B_{n,k}$ 
is a  number of particles in the
$\Lambda$-coalescent started with $n$ particles right after the number of
particles drops below to~$k+1$.  
Then for any $k\le n$
the total number of collisions is decomposable as
\[
C_n=_d C_{n\fromto k}+C_{B_{n,k}}^{(1)}
\]
where in the RHS  $(C_b^{(1)})$ is an independent copy of $(C_b)$.  
This can be iterated as
\begin{equation}\label{eq:cdisteq2}
C_n= C_{n\fromto k_1}+C_{B_{n,k_1}\fromto k_2}^{(1)}+C_{B_{n,k_2}\fromto
k_3}^{(2)}+\dots+C_{B_{n,k_{\ell-1}}\fromto k_\ell}^{(\ell)}
\end{equation}
for any finite sequence $k_\ell\le k_{\ell-1}\le\dots\le k_1\le n$, with
the convention that $C_{b\fromto k}=0$ 
for $b\le k$.

Suppose that Lemma~\ref{lem:Cntondelta} holds for some $\upsilon<1$.
Let $\ell=\intfloor{\frac{-\log(2-\alpha)}{\log\upsilon}}$, then 
$\upsilon^{\ell+1}<1/(2-\alpha)$.

For each $m=1,\dots,\ell$, by Lemma~\ref{lem:Cntondelta} applied with
$\intfloor{n^{\upsilon^m}}$ instead of $n$, 
\[
\prob\left[\frac{C_{\intfloorsm{n^{\upsilon^m}}\fromto
\intfloorsm{n^{\upsilon^{m+1}}}}
 -\left(\intfloorsm{n^{\upsilon^m}}-\intfloorsm{n^{\upsilon^{m+1}}}
\right)(1-\alpha) }
{n^{{\upsilon^m}/(2-\alpha)}(1-\alpha)}\le
t\right]\to F_{2-\alpha}(t).
\]
Consequently, since $\upsilon<1$, for all $m>0$ 
\[
\frac{C_{\intfloorsm{n^{\upsilon^m}}\fromto
\intfloorsm{n^{\upsilon^{m+1}}}}-\left(\intfloorsm{n^{\upsilon^m}}
-\intfloorsm{n^{\upsilon^{m+1}}}
\right)(1-\alpha)}
 {n^{1/(2-\alpha)}(1-\alpha)} \to_d \bfdelta_0,\qquad
n\to\infty,
\]
where $\bfdelta_0$ is the $\delta$-measure at zero.  Moreover,
for some fixed $\tau\in{]0,1/(2-\alpha)[}$ starting the coalescent with
$\intfloorsm{n^{\upsilon^k}-n^\tau}$ particles 
instead of $\intfloorsm{n^{\upsilon^m}}$ particles results in 
the same asymptotic behaviour:
\[
\frac{C_{\intfloorsm{n^{\upsilon^m}-n^\tau}\fromto
\intfloorsm{n^{\upsilon^{m+1}}}}-\left(\intfloorsm{n^{\upsilon^m}}
-\intfloorsm{n^{\upsilon^{m+1}}}
\right)(1-\alpha)}
 {n^{1/(2-\alpha)}(1-\alpha)} \to_d \bfdelta_0,\qquad n\to\infty.
\]
Denote by $E_m$ the event $B_{n,\intfloorsm{n^{\upsilon^m}}}\ge
n^{\upsilon^m}-n^\tau$, i.e.\ that the Markov process
$\CM_n$ undershoots $n^{\upsilon^m}$ not more than by $n^\tau$. 
Thus, given $E_m$, 
\[
\frac{C_{B_{n,n^{\upsilon^m}}\fromto
\intfloorsm{n^{\upsilon^{m+1}}}}-\left(\intfloorsm{n^{\upsilon^m}}
-\intfloorsm{n^{\upsilon^{m+1}}}
\right)(1-\alpha)}
 {n^{1/(2-\alpha)}(1-\alpha)} \to_d \bfdelta_0,\qquad n\to\infty.
\]
by monotonicity of the number of collisions.  
Using \eqref{eq:cdisteq2} with
$k_m=\intfloorsm{n^{\upsilon^m}}$  and Lemma~\ref{lem:Cntondelta} yields
the desired convergence for $C_n$, 
given $E_m$ holds for all $m=1,\dots,\ell$.

Lemma~\ref{lem:latail} ensures that for any $\tau>0$ the 
 probability of $E_m$ grows to~1, as
$n\to\infty$, for all $m=1,\dots,\ell$.  Hence
$\prob\left[\cap_{m=1}^\ell E_m\right]\to1$ and the convergence in
distribution conditioned on
$\cap_{m=1}^\ell E_m$ is equivalent to the unconditional convergence in
distribution. 
\end{proof}

\end{document}